\documentclass{article}%
\usepackage{amsfonts}
\usepackage{amssymb}
\usepackage{amsmath}%
\usepackage{graphicx}

\numberwithin{equation}{section}
\begin{document}

\title{An involution for recurrences of Tricomi and other special functions}
\author{John T. Conway\thanks{CONTACT John T. Conway Email: john.conway@uia.no}\\Department of Engineering and Science,\\University of Agder, Grimstad, Norway}
\maketitle
\date{}

\begin{abstract}
Recently two new integrated recurrences were given for the Kummer function
${}_{1}F_{1}\left(  a;b;x\right)  $ which together with the six recurrences
already given in standard literature gave a total of eight recurrences, each
one linking the baseline function ${}_{1}F_{1}\left(  a;b;x\right)  $ to one
of it's eight contigous functions. Two new recurrences for the Tricomi
function $U\left(  a;b;x\right)  $ are given here which complete the set of
eight recurrences for $U\left(  a;b;x\right)  $\ linking it to it's contiguous
functions. The involution presented previously is examined in detail for the
recurrences of Cylinder, Associated Legendre and Kummer/Tricomi functions.
This involution converts any given recurrence to a recurrence with the
parameter shifts reversed, which creates conjugate pairs. One pair for
Cylinder functions, two pairs for Associated Legendre functions and 4 pairs
each for Kummer and Tricomi functions. Uniqueness is proven using Wronskians
for the integrated recurrence relations, such that contiguous functions are
single occupancy slots with only one integrated recurrence possible for each
slot and for each independent solution of the ODE. For Cylinder and Associated
Legendre functions this is simple, as the recurrences are identical for each
independent solution. However, for Kummer and Tricomi functions the
recurrences are different. This difficulty is bypassed using Gamma functions
to scale the Tricomi recurrences to match the Kummer recurrences for each
conjugate pair. It is not possible to match globally. This allows uniqueness
to be proven using Wronskians on a pair by pair basis. Hence the involution
system behaves like a Pauli Exclusion principle for contiguous functions, with
only one recurrence possible for each solution of the ODE.

\begin{description}
\item[Keywords:] Differential equations; Kummer functions; Tricomi functions;
Associated Legendre functions; Bessel functions;

\item[AMS\ Subject Classification:] 34B30; 33C10; 33C15; 42C10; 33C45;

\end{description}
\end{abstract}

\section{Introduction}

Many special functions $y\left(  x\right)  $ obey a second-order linear
homogeneous differential equation of the form%
\begin{equation}
y^{\prime\prime}\left(  x\right)  +p\left(  x\right)  y^{\prime}\left(
x\right)  +q\left(  x\right)  y\left(  x\right)  =0\label{Eqn1}%
\end{equation}
and these functions usually possess well known three term recurrence relations
for the parameters used in the definition of the functions. For example, the
Bessel equation for the general cylinder function $Z_{n}\left(  x\right)
\equiv A_{1}J_{n}\left(  x\right)  +A_{2}Y_{n}\left(  x\right)  $ is%
\begin{equation}
Z_{n}^{\prime\prime}\left(  x\right)  +\frac{1}{x}Z_{n}^{\prime}\left(
x\right)  +\left(  1-\frac{n^{2}}{x^{2}}\right)  Z_{n}\left(  x\right)
\label{Eqn2}%
\end{equation}
and in this case there is only one parameter $n$. The the well-known three
term recurrences are
\begin{equation}
Z_{n}^{\prime}\left(  x\right)  +\frac{n}{x}Z_{n}\left(  x\right)
=Z_{n-1}\left(  x\right)  \label{Eqn3}%
\end{equation}%
\begin{equation}
Z_{n}^{\prime}\left(  x\right)  -\frac{n}{x}Z_{n}\left(  x\right)
=-Z_{n+1}\left(  x\right)  \label{Eqn4}%
\end{equation}
and these relations can be transformed with integrating factors to give the
two-term integrated recurrences%
\begin{equation}
\left(  x^{n}Z_{n}\left(  x\right)  \right)  ^{\prime}=x^{n}Z_{n-1}\left(
x\right)  \label{Eqn5}%
\end{equation}%
\begin{equation}
\left(  x^{-n}Z_{n}\left(  x\right)  \right)  ^{\prime}=-x^{-n}Z_{n+1}\left(
x\right)  \text{.}\label{Eqn6}%
\end{equation}
For Equation (\ref{Eqn1}) the differential equation for the function $u\left(
x\right)  =c\left(  x\right)  y\left(  x\right)  $ is given by the general
transformation equation [1,2,3]%
\[
u^{\prime\prime}\left(  x\right)  +\left(  p\left(  x\right)  -2\frac
{c^{\prime}\left(  x\right)  }{c\left(  x\right)  }\right)  u^{\prime}\left(
x\right)
\]%
\begin{equation}
+\left(  \left(  \frac{c^{\prime}\left(  x\right)  }{c\left(  x\right)
}\right)  ^{2}-\left(  \frac{c^{\prime}\left(  x\right)  }{c\left(  x\right)
}\right)  ^{\prime}-p\left(  x\right)  \frac{c^{\prime}\left(  x\right)
}{c\left(  x\right)  }+q\left(  x\right)  \right)  u\left(  x\right)
=0\text{.}\label{Eqn7}%
\end{equation}
The left-hand side of either integrated recurrence can be written in the form
$u_{n}\left(  x\right)  \equiv\left(  c_{n}\left(  x\right)  Z_{n}\left(
x\right)  \right)  ^{\prime}$ and the factor $c_{n}\left(  x\right)  $ can be
substituted into Equation (\ref{Eqn7}) to give differential equations for
$u_{n}\left(  x\right)  $. For equations (\ref{Eqn5}) and (\ref{Eqn6}) this
gives the well known equations
\begin{equation}
u_{n}^{\prime\prime}\left(  x\right)  +\frac{1\mp2n}{x}u_{n}^{\prime}\left(
x\right)  +u_{n}\left(  x\right)  =0\label{Eqn8}%
\end{equation}
with the general solutions
\begin{equation}
u_{n}\left(  x\right)  =x^{\pm n}Z_{n}\left(  x\right)  .\label{Eqn9}%
\end{equation}

\subsection{The Involution}

The integrating factor for the two leftmost terms of Equation (\ref{Eqn8}) is
$f\left(  x\right)  =x^{1\mp2n}$ and applying this to the differential
equation gives%
\begin{equation}
\left[  x^{1\mp2n}u_{n}^{\prime}\left(  x\right)  \right]  ^{\prime}%
=-x^{1\mp2n}u_{n}\left(  x\right) \label{Eqn10}%
\end{equation}
and from Equation (\ref{Eqn9}) this becomes%
\begin{equation}
\left[  x^{1\mp2n}\left(  x^{\pm n}Z_{n}\left(  x\right)  \right)  ^{\prime
}\right]  ^{\prime}=-x^{1\mp n}Z_{n}\left(  x\right)  \text{.}\label{Eqn11}%
\end{equation}
Choosing $u_{n}\left(  x\right)  =x^{n}Z_{n}\left(  x\right)  $ in Equation
(\ref{Eqn11}) gives%
\begin{equation}
\left[  x^{1-2n}\left(  x^{n}Z_{n}\left(  x\right)  \right)  ^{\prime}\right]
^{\prime}=-x^{1-n}Z_{n}\left(  x\right) \label{Eqn12}%
\end{equation}
and substituting the recurrence (\ref{Eqn5}) into the left-hand side of
Equation (\ref{Eqn12}) integrates this equation once to give%
\begin{equation}
\left[  x^{1-n}Z_{n-1}\left(  x\right)  \right]  ^{\prime}=-x^{1-n}%
Z_{n}\left(  x\right)  \text{.}\label{Eqn13}%
\end{equation}
Shifting the parameter in Equation (\ref{Eqn13}) such that $n\rightarrow n+1$
gives%
\begin{equation}
\left[  x^{-n}Z_{n}\left(  x\right)  \right]  ^{\prime}=-x^{-n}Z_{n+1}\left(
x\right) \label{Eqn14}%
\end{equation}
which is identical to Equation (\ref{Eqn6}). Similarly, choosing $u_{n}\left(
x\right)  =x^{-n}Z_{n}\left(  x\right)  $ in Equation (\ref{Eqn11}) and
following the same procedure leads to%
\begin{equation}
\left[  x^{n}Z_{n}\left(  x\right)  \right]  ^{\prime}=x^{n}Z_{n-1}\left(
x\right) \label{Eqn15}%
\end{equation}
which is identical to Equation (\ref{Eqn5}).

Employing each recurrence in Equation (\ref{Eqn11}) necessarily gives the
other recurrence because the three-term recurrences (\ref{Eqn3}) and
(\ref{Eqn4}) are both unique for their particular parameter shifts and hence
the integrated recurrences in Equations (\ref{Eqn5}) and (\ref{Eqn6}) are also
unique. This can be shown very easily using Wronskians.\newline

\subsubsection{Proof of Uniqueness}

Let us suppose the existence of two hypothetical recurrences of the same form
as Equation (\ref{Eqn3}), namely%
\begin{equation}
Z_{n}^{\prime}\left(  x\right)  +g\left(  x\right)  Z_{n}\left(  x\right)
=h\left(  x\right)  Z_{n-1}\left(  x\right) \label{Eqn16}%
\end{equation}%
\begin{equation}
Z_{n}^{\prime}\left(  x\right)  +r\left(  x\right)  Z_{n}\left(  x\right)
=s\left(  x\right)  Z_{n-1}\left(  x\right) \label{Eqn17}%
\end{equation}
Then the terms in $Z_{n-1}\left(  x\right)  $ can be eliminated between the
two recurrences to give%
\begin{equation}
\frac{Z_{n}^{\prime}\left(  x\right)  }{Z_{n}\left(  x\right)  }=\alpha\left(
x\right) \label{Eqn18}%
\end{equation}
where%
\begin{equation}
\alpha\left(  x\right)  =\frac{h\left(  x\right)  r\left(  x\right)  -s\left(
x\right)  g\left(  x\right)  }{s\left(  x\right)  -h\left(  x\right)
}.\label{Eqn19}%
\end{equation}
Since the function $\alpha\left(  x\right)  $ is the same for both independent
solutions of Equation (\ref{Eqn2}), then Equation (\ref{Eqn18}) gives
immediately%
\[
\frac{J_{n}^{\prime}(x)}{J_{n}\left(  x\right)  }=\frac{Y_{n}^{\prime}\left(
x\right)  }{Y_{n}(x)}%
\]
and this gives the Wronskian $W\left(  J_{n}\left(  x\right)  ,Y_{n}\left(
x\right)  \right)  \equiv0$, which is a contradiction. Hence Equations
(\ref{Eqn3}) and (\ref{Eqn4}) are the only possible three term recurrences for
their respective parameter shifts. Consequently, Equations (\ref{Eqn5}) and
(\ref{Eqn6}) are the only corresponding integrated recurrences and form a
conjugate pair under the action of the involution. A key element of this
particular proof is that the recurrences must be obeyed by both solutions of
the ODE. This is also true for the Associated Legendre functions treated
below, but is not true for the Kummer and Tricomi functions treated later,
where the proof has to be modified.

\subsection{Associated Legendre functions}

The general Associated Legendre function $y\left(  x\right)  =R_{\nu}^{\mu
}\left(  x\right)  \equiv B_{1}P_{\nu}^{\mu}\left(  x\right)  +B_{2}Q_{\nu
}^{\mu}\left(  x\right)  $ obeys the equation%
\begin{equation}
y^{\prime\prime}\left(  x\right)  -\frac{2x}{1-x^{2}}y^{\prime}\left(
x\right)  +\left(  \frac{\nu\left(  \nu+1\right)  }{1-x^{2}}-\frac{\mu^{2}%
}{\left(  1-x^{2}\right)  ^{2}}\right)  y\left(  x\right)  =0.\label{Eqn20}%
\end{equation}
As in the case of cylinder functions, the four standard recurrence relations
are obeyed by both independent solutions of Equation (\ref{Eqn20}), and hence
uniqueness of the three-term recurrences can be proven as above using
Wronskians. The integrated recurrences, transformed equations and the
integrating factors for the transformed equations were given in [4] in a quite
different context and for convenience these results are repeated below.%
\begin{equation}
\left[  \left(  1-x^{2}\right)  ^{\left(  \nu+1\right)  /2}R_{\nu}^{\mu
}\left(  x\right)  \right]  ^{\prime}=-\left(  \nu-\mu+1\right)  \left(
1-x^{2}\right)  ^{\left(  \nu-1\right)  /2}R_{\nu+1}^{\mu}\left(  x\right)
\label{Eqn21}%
\end{equation}%
\begin{equation}
\left[  \left(  1-x^{2}\right)  ^{-\nu/2}R_{\nu}^{\mu}\left(  x\right)
\right]  ^{\prime}=\left(  \nu+\mu\right)  \left(  1-x^{2}\right)  ^{-\nu
/2-1}R_{\nu-1}^{\mu}\left(  x\right) \label{Eqn22}%
\end{equation}%
\begin{equation}
\left[  \left(  1-x^{2}\right)  ^{-\mu/2}R_{\nu}^{\mu}\left(  x\right)
\right]  ^{\prime}=-\left(  1-x^{2}\right)  ^{-\left(  \mu+1\right)  /2}%
R_{\nu}^{\mu+1}\left(  x\right) \label{Eqn23}%
\end{equation}%
\begin{equation}
\left[  \left(  1-x^{2}\right)  ^{\mu/2}R_{\nu}^{\mu}\left(  x\right)
\right]  ^{\prime}=\left(  \nu-\mu+1\right)  \left(  \nu+\mu\right)  \left(
1-x^{2}\right)  ^{\left(  \mu-1\right)  /2}R_{\nu}^{\mu-1}\left(  x\right)
\text{.}\label{Eqn24}%
\end{equation}
We can define%
\begin{equation}
z_{1}\left(  x\right)  =\left(  1-x^{2}\right)  ^{\left(  \nu+1\right)
/2}R_{\nu}^{\mu}\left(  x\right) \label{Eqn25}%
\end{equation}%
\begin{equation}
z_{2}=\left(  1-x^{2}\right)  ^{-\nu/2}R_{\nu}^{\mu}\left(  x\right)
\label{Eqn26}%
\end{equation}%
\begin{equation}
z_{3}\left(  x\right)  =\left(  1-x^{2}\right)  ^{-\mu/2}R_{\nu}^{\mu}\left(
x\right) \label{Eqn27}%
\end{equation}%
\begin{equation}
z_{4}\left(  x\right)  =\left(  1-x^{2}\right)  ^{\mu/2}R_{\nu}^{\mu}\left(
x\right) \label{Eqn28}%
\end{equation}
The transformed equations for the $z_{i}\left(  x\right)  $ are given by%
\begin{equation}
z_{1}^{\prime\prime}\left(  x\right)  +\frac{2\nu x}{1-x^{2}}z_{1}^{\prime
}\left(  x\right)  =-\frac{\left(  \nu+1\right)  ^{2}-\mu^{2}}{\left(
1-x^{2}\right)  ^{2}}z_{1}\left(  x\right) \label{Eqn29}%
\end{equation}%
\begin{equation}
z_{2}^{\prime\prime}\left(  x\right)  -\frac{2\left(  \nu+1\right)  x}%
{1-x^{2}}z_{2}^{\prime}\left(  x\right)  =-\frac{\nu^{2}-\mu^{2}}{\left(
1-x^{2}\right)  ^{2}}z_{2}\left(  x\right) \label{Eqn30}%
\end{equation}%
\begin{equation}
z_{3}^{\prime\prime}\left(  x\right)  -\frac{2\left(  1+\mu\right)  x}%
{1-x^{2}}z_{3}^{\prime}\left(  x\right)  =-\frac{\nu\left(  \nu+1\right)
-\mu\left(  \mu+1\right)  }{1-x^{2}}z_{3}\left(  x\right) \label{Eqn31}%
\end{equation}%
\begin{equation}
z_{4}^{\prime\prime}\left(  x\right)  -\frac{2\left(  1-\mu\right)  x}%
{1-x^{2}}z_{4}^{\prime}\left(  x\right)  =-\frac{\nu\left(  \nu+1\right)
-\mu\left(  \mu-1\right)  }{1-x^{2}}z_{4}\left(  x\right) \label{Eqn32}%
\end{equation}
and the integrating factors for the left-hand sides of these equations are,
respectively%
\begin{equation}
f_{1}\left(  x\right)  =\exp\left(
{\displaystyle\int}
\frac{2\nu x}{1-x^{2}}dx\right)  =\left(  1-x^{2}\right)  ^{-\nu}\label{Eqn33}%
\end{equation}%
\begin{equation}
f_{2}\left(  x\right)  =\exp\left(  -%
{\displaystyle\int}
\frac{2\left(  \nu+1\right)  x}{1-x^{2}}dx\right)  =\left(  1-x^{2}\right)
^{\nu+1}\label{Eqn34}%
\end{equation}%
\begin{equation}
f_{3}\left(  x\right)  =\exp\left(  -%
{\displaystyle\int}
\frac{2\left(  1+\mu\right)  x}{1-x^{2}}dx\right)  =\left(  1-x^{2}\right)
^{1+\mu}\label{Eqn35}%
\end{equation}%
\begin{equation}
f_{4}\left(  x\right)  =\exp\left(  -%
{\displaystyle\int}
\frac{2\left(  1-\mu\right)  x}{1-x^{2}}dx\right)  =\left(  1-x^{2}\right)
^{1-\mu}\label{Eqn36}%
\end{equation}

Applying Equations (\ref{Eqn25}) and (\ref{Eqn33}) in Equation (\ref{Eqn29})
gives%
\begin{equation}
\left[  \left(  1-x^{2}\right)  ^{-\nu}\left(  \left(  1-x^{2}\right)
^{\left(  \nu+1\right)  /2}R_{\nu}^{\mu}\left(  x\right)  \right)  ^{\prime
}\right]  ^{\prime}=-\left(  \left(  \nu+1\right)  ^{2}-\mu^{2}\right)
\left(  1-x^{2}\right)  ^{\left(  \nu-3\right)  /2}R_{\nu}^{\mu}\left(
x\right) \label{Eqn37}%
\end{equation}
and employing Equation (\ref{Eqn21}) in Equation (\ref{Eqn37}) gives after
some reduction%
\begin{equation}
\left[  \left(  1-x^{2}\right)  ^{-\left(  \nu+1\right)  /2}R_{\nu+1}^{\mu
}\left(  x\right)  \right]  ^{\prime}=\left(  \nu+\mu+1\right)  \left(
1-x^{2}\right)  ^{-\left(  \nu+3\right)  /2}R_{\nu}^{\mu}\left(  x\right)
\text{.}\label{Eqn38}%
\end{equation}
Shifting the parameter $\nu$ such that $\nu\rightarrow\nu-1$ then gives
\begin{equation}
\left[  \left(  1-x^{2}\right)  ^{-\nu/2}R_{\nu}^{\mu}\left(  x\right)
\right]  ^{\prime}=\left(  \nu+\mu\right)  \left(  1-x^{2}\right)  ^{-\nu
/2-1}R_{\nu-1}^{\mu}\left(  x\right) \label{Eqn39}%
\end{equation}
which is identical to Equation (\ref{Eqn22}). Similarly, applying Equations
(\ref{Eqn26}) and (\ref{Eqn34}) in Equation (\ref{Eqn30}) gives%
\begin{equation}
\left[  \left(  1-x^{2}\right)  ^{\nu+1}\left(  \left(  1-x^{2}\right)
^{-\nu/2}R_{\nu}^{\mu}\left(  x\right)  \right)  ^{\prime}\right]  ^{\prime
}=-\left(  \nu^{2}-\mu^{2}\right)  \left(  1-x^{2}\right)  ^{\nu/2-1}R_{\nu
}^{\mu}\left(  x\right) \label{Eqn40}%
\end{equation}
and employing Equation (\ref{Eqn22}) in this equation gives%
\begin{equation}
\left[  \left(  1-x^{2}\right)  ^{\nu/2}R_{\nu-1}^{\mu}\left(  x\right)
\right]  ^{\prime}=-\left(  \nu-\mu\right)  \left(  1-x^{2}\right)  ^{\nu
/2-1}R_{\nu}^{\mu}\left(  x\right) \label{Eqn41}%
\end{equation}
and shifting the parameter $\nu$ such that $\nu\rightarrow\nu+1$ in Equation
(\ref{Eqn41}) gives%
\begin{equation}
\left[  \left(  1-x^{2}\right)  ^{\left(  \nu+1\right)  /2}R_{\nu}^{\mu
}\left(  x\right)  \right]  ^{\prime}=-\left(  \nu-\mu+1\right)  \left(
1-x^{2}\right)  ^{\left(  \nu-1\right)  /2}R_{\nu+1}^{\mu}\left(  x\right)
\label{Eqn42}%
\end{equation}
which is identical to Equation (\ref{Eqn21}). Hence Equations (\ref{Eqn21}%
)\ and (\ref{Eqn22}) form a conjugate pair for the parameter $\nu$.

Applying Equations (\ref{Eqn27}) and (\ref{Eqn35}) in Equation (\ref{Eqn31})
gives%
\begin{equation}
\left[  \left(  1-x^{2}\right)  ^{1+\mu}\left(  \left(  1-x^{2}\right)
^{-\mu/2}R_{\nu}^{\mu}\left(  x\right)  \right)  ^{\prime}\right]  ^{\prime
}=-\left(  \nu\left(  \nu+1\right)  -\mu\left(  \mu+1\right)  \right)  \left(
1-x^{2}\right)  ^{\mu/2}R_{\nu}^{\mu}\left(  x\right) \label{Eqn43}%
\end{equation}
and applying Equation (\ref{Eqn23}) in this equation gives%
\begin{equation}
\left[  \left(  1-x^{2}\right)  ^{\left(  \mu+1\right)  /2}R_{\nu}^{\mu
+1}\left(  x\right)  \right]  ^{\prime}=\left(  \nu\left(  \nu+1\right)
-\mu\left(  \mu+1\right)  \right)  \left(  1-x^{2}\right)  ^{\mu/2}R_{\nu
}^{\mu}\left(  x\right)  \text{.}\label{Eqn44}%
\end{equation}
Shifting the parameter $\mu$ such that $\mu\rightarrow\mu-1$ gives%
\begin{equation}
\left[  \left(  1-x^{2}\right)  ^{\mu/2}R_{\nu}^{\mu}\left(  x\right)
\right]  ^{\prime}=\left(  \nu-\mu+1\right)  \left(  \nu+\mu\right)  \left(
1-x^{2}\right)  ^{\left(  \mu-1\right)  /2}R_{\nu}^{\mu-1}\left(  x\right)
\label{Eqn45}%
\end{equation}
which is identical to Equation (\ref{Eqn24}). Similarly, applying Equations
(\ref{Eqn28}) and (\ref{Eqn36}) in Equation (\ref{Eqn32}) gives%
\begin{equation}
\left[  \left(  1-x^{2}\right)  ^{1-\mu}\left(  \left(  1-x^{2}\right)
^{\mu/2}R_{\nu}^{\mu}\left(  x\right)  \right)  ^{\prime}\right]  ^{\prime
}=-\left(  \nu\left(  \nu+1\right)  -\mu\left(  \mu-1\right)  \right)  \left(
1-x^{2}\right)  ^{-\mu/2}R_{\nu}^{\mu}\left(  x\right) \label{Eqn46}%
\end{equation}
and employing Equation (\ref{Eqn24}) in Equation (\ref{Eqn46}) gives, after
cancelling the constants%
\begin{equation}
\left[  \left(  1-x^{2}\right)  ^{\left(  1-\mu\right)  /2}R_{\nu}^{\mu
-1}\left(  x\right)  \right]  ^{\prime}=-\left(  1-x^{2}\right)  ^{-\mu
/2}R_{\nu}^{\mu}\left(  x\right)  \text{.}\label{Eqn47}%
\end{equation}
Finally, shifting the parameter $\mu$ such that $\mu\rightarrow\mu+1$ gives%
\begin{equation}
\left[  \left(  1-x^{2}\right)  ^{-\mu/2}R_{\nu}^{\mu}\left(  x\right)
\right]  ^{\prime}=-\left(  1-x^{2}\right)  ^{-\left(  \mu+1\right)  /2}%
R_{\nu}^{\mu+1}\left(  x\right) \label{Eqn48}%
\end{equation}
which is identical to Equation (\ref{Eqn31}). Hence Equations (\ref{Eqn45})
and (\ref{Eqn48}) form a second conjugate pair for shifts in the parameter
$\mu$.

The next section applies the framework above to Tricomi and Kummer functions,
which are two independent solutions of Kummer's equation. These two functions
have different recurrence relations and must in most ways be treated
separately. This gives an initial difficulty for the Wronskian proof of
uniqueness given above but this difficulty can be bypassed in a very nice and
interesting manner. Results have already been presented in [3] for the Kummer
function $_{1}F_{1}\left(  a;b;x\right)  $ which included deriving two new
integrated recurrences for $_{1}F_{1}\left(  a;b;x\right)  $ to make a total
of 8, one for each of the contiguous functions. The first of these new
recurrences was derived from some contiguous relations given in [5] and the
other was obtained by applying the involution to this relation. The next
section will derive the corresponding results for the Tricomi function
$U\left(  a;b;x\right)  $ and then combine results for both functions to prove uniqueness.

\section{Kummer and Tricomi functions}

The Kummer equation is%
\begin{equation}
y^{\prime\prime}\left(  x\right)  +\left(  \frac{b}{x}-1\right)  y^{\prime
}\left(  x\right)  -\frac{a}{x}y\left(  x\right)  =0\label{EqnA1}%
\end{equation}
and two independent solutions of this equation are the Kummer confluent
hypergeometric function $_{1}F_{1}\left(  a;b;x\right)  $ and the Tricomi
confluent hypergeometric function $U\left(  a;b;x\right)  $. The primary
differential relation for the Tricomi function is [6]:
\begin{equation}
\left[  U\left(  a;b;x\right)  \right]  ^{\prime}=-aU\left(  a+1;b+1;x\right)
\label{EqnA2}%
\end{equation}

\subsection{Tricomi integrated recurrences}

Six three-term contiguous relations are given in [6] and four of these are:%
\begin{equation}
U\left(  a-1;b;x\right)  +\left(  b-2a-x\right)  U\left(  a;b;x\right)
+a\left(  a-b+1\right)  U\left(  a+1;b;x\right)  =0\label{EqnA3}%
\end{equation}%
\begin{equation}
\left(  b-a-1\right)  U\left(  a;b-1;x\right)  +\left(  1-b-x\right)  U\left(
a;b;x\right)  +xU\left(  a;b+1;x\right)  =0\label{EqnA4}%
\end{equation}%
\begin{equation}
U\left(  a;b;x\right)  -aU\left(  a+1;b;x\right)  -U\left(  a;b-1;x\right)
=0\label{EqnA5}%
\end{equation}%
\begin{equation}
\left(  a-1+x\right)  U\left(  a;b;x\right)  -U\left(  a-1;b;x\right)
+\left(  a-b+1\right)  U\left(  a;b-1;x\right) \label{EqnA6}%
\end{equation}
Shifting the parameter $b$ such that $b\rightarrow b+1$ in Equations
(\ref{EqnA3}) and (\ref{EqnA5}) gives, respectively%
\begin{equation}
U\left(  a-1;b+1;x\right)  +\left(  b-2a+1-x\right)  U\left(  a;b+1;x\right)
+a\left(  a-b\right)  U\left(  a+1;b+1;x\right)  =0\label{EqnA7}%
\end{equation}%
\begin{equation}
U\left(  a;b+1;x\right)  -aU\left(  a+1;b+1;x\right)  -U\left(  a;b;x\right)
=0\label{EqnA8}%
\end{equation}
Combining these two equations to eliminate the terms in $U\left(
a;b+1;x\right)  $ gives the contiguous relation
\begin{equation}
U\left(  a-1;b+1;x\right)  -a\left(  a-1+x\right)  U\left(  a+1;b+1;x\right)
+\left(  b-2a+1-x\right)  U\left(  a;b;x\right)  =0\label{EqnA9}%
\end{equation}
and employing Equation (\ref{EqnA2}) in Equation (\ref{EqnA9}) gives%
\begin{equation}
\left[  U\left(  a;b;x\right)  \right]  ^{\prime}+\frac{\left(
b-2a+1-x\right)  }{\left(  a-1+x\right)  }U\left(  a;b;x\right)
+\frac{U\left(  a-1;b+1;x\right)  }{\left(  a-1+x\right)  }=0\label{EqnA10}%
\end{equation}
The integrating factor for this differential recurrence is%
\begin{equation}
\exp\left(
{\displaystyle\int}
\frac{\left(  b-2a+1-x\right)  }{\left(  a-1+x\right)  }dx\right)  =\left(
a-1+x\right)  ^{b-a}e^{-x}\label{EqnA11}%
\end{equation}
so a new integrated recurrence for the contiguous function $U\left(
a-1;b+1;x\right)  $ is given by
\begin{equation}
\left[  \left(  a-1+x\right)  ^{b-a}e^{-x}U\left(  a;b;x\right)  \right]
^{\prime}=-\left(  a-1+x\right)  ^{b-a-1}e^{-x}U\left(  a-1;b+1;x\right)
\text{.}\label{EqnA12}%
\end{equation}

Similarly, we can shift the parameter $a$ such that $a\rightarrow a+1$ in
Equations (\ref{EqnA4}) and (\ref{EqnA5}) which gives, respectively%
\begin{equation}
\left(  b-a-2\right)  U\left(  a+1;b-1;x\right)  +\left(  1-b-x\right)
U\left(  a+1;b;x\right)  +xU\left(  a+1;b+1;x\right)  =0\label{EqnA13}%
\end{equation}%
\begin{equation}
\left(  a+x\right)  U\left(  a+1;b;x\right)  -U\left(  a;b;x\right)  +\left(
a-b+2\right)  U\left(  a+1;b-1;x\right)  =0\label{EqnA14}%
\end{equation}
Multiplying Equation (\ref{EqnA13}) by $a+x$ and Equation (\ref{Eqn14}) by
$1-b-x$ to eliminate the term in $U\left(  a+1,b,x\right)  $ gives the
recurrence%
\[
\left(  a-b+1\right)  (b-a-2)U\left(  a+1;b-1;x\right)  +\left(  1-b-x\right)
U\left(  a;b;x\right)
\]%
\begin{equation}
+x\left(  a+x\right)  U\left(  a+1;b+1;x=0\right)  \label{EqnA15}%
\end{equation}
Employing the derivative relation (\ref{EqnA2}) in (\ref{EqnA15}) gives the
differential recurrence%
\begin{equation}
U^{\prime}\left(  a;b;x\right)  +\frac{a\left(  x+b-1\right)  }{x\left(
a+x\right)  }U\left(  a;b;x\right)  =\frac{a\left(  a-b+1\right)
(b-a-2)}{x\left(  a+x\right)  }U\left(  a+1;b-1;x\right)  \label{EqnA16}%
\end{equation}
and the integrating factor for this equation is%
\begin{equation}
\exp\left(
{\displaystyle\int}
\frac{a\left(  x+b-1\right)  }{x\left(  a+x\right)  }dx\right)  =x^{b-1}%
\left(  a+x\right)  ^{a-b+1}\text{.}\label{EqnA17}%
\end{equation}
Hence Equation (\ref{EqnA16}) can be integrated to give the integrated
recurrence%
\begin{equation}
\left[  x^{b-1}\left(  a+x\right)  ^{a-b+1}U\left(  a;b;x\right)  \right]
^{\prime}=a(b-a-2)\left(  a-b+1\right)  x^{b-2}\left(  a+x\right)
^{a-b}U\left(  a+1;b-1;x\right)  \text{.}\label{EqnA18}%
\end{equation}
Equations (\ref{EqnA12}) and (\ref{EqnA18}) are new integrated recurrences for
the contiguous functions $U\left(  a-1;b+1;x\right)  $ and $U\left(
a+1;b-1;x\right)  $. The derivations of these recurrences given above from
well known recurrences means they have been hiding in plain sight for a long time.

There are five additional integrated differential recurrences given in [6,7,8]
and combining these with Equations (\ref{EqnA2}), (\ref{EqnA12}) and
(\ref{EqnA18}) gives the full set of eight integrated recurrence relations for
the Tricomi functions%

\begin{equation}
\left[  U\left(  a;b;x\right)  \right]  ^{\prime}=-aU\left(  a+1;b+1;x\right)
\label{EqnA19}%
\end{equation}%
\begin{equation}
\left[  x^{a}U\left(  a;b;x\right)  \right]  ^{\prime}=a\left(  a-b+1\right)
x^{a-1}U\left(  a+1;b;x\right)  \label{EqnA20}%
\end{equation}%
\begin{equation}
\left[  x^{b-1}U\left(  a;b;x\right)  \right]  ^{\prime}=-\left(
a-b+1\right)  x^{b-2}U\left(  a;b-1;x\right)  \label{EqnA21}%
\end{equation}%
\begin{equation}
\left[  x^{b-a}e^{-x}U\left(  a;b;x\right)  \right]  ^{\prime}=-x^{b-a-1}%
e^{-x}U\left(  a-1;b;x\right)  \label{EqnA22}%
\end{equation}%
\begin{equation}
\left[  e^{-x}U\left(  a;b;x\right)  \right]  ^{\prime}=-e^{-x}U\left(
a;b+1;x\right)  \label{EqnA23}%
\end{equation}%
\begin{equation}
\left[  x^{b-1}e^{-x}U\left(  a;b;x\right)  \right]  ^{\prime}=-x^{b-2}%
e^{-x}U\left(  a-1;b-1;x\right)  \label{EqnA24}%
\end{equation}%
\begin{equation}
\left[  \left(  a-1+x\right)  ^{b-a}e^{-x}U\left(  a;b;x\right)  \right]
^{\prime}=-\left(  a-1+x\right)  ^{b-a-1}e^{-x}U\left(  a-1;b+1;x\right)
\label{EqnA25}%
\end{equation}%
\[
\left[  x^{b-1}\left(  a+x\right)  ^{a-b+1}U\left(  a;b;x\right)  \right]
^{\prime}=
\]%
\begin{equation}
a(b-a-2)\left(  a-b+1\right)  x^{b-2}\left(  a+x\right)  ^{a-b}U\left(
a+1;b-1;x\right)  \text{.}\label{EqnA26}%
\end{equation}

The eight corresponding integrated recurrences for Kummer function $_{1}%
F_{1}\left(  a;b;x\right)  $ were given in [3] and are repeated below:
\begin{equation}
\left[  {}_{1}F_{1}\left(  a;b;x\right)  \right]  ^{\prime}=\frac{a}%
{b}\text{\smallskip}_{1}F_{1}\left(  a+1;b+1;x\right)  \label{EqnA27}%
\end{equation}%
\begin{equation}
\left[  x^{a}{}_{1}F_{1}\left(  a;b;x\right)  \right]  ^{\prime}=ax^{a-1}%
{}_{1}F_{1}\left(  a+1;b;x\right)  \label{EqnA28}%
\end{equation}%
\begin{equation}
\left[  x^{b-1}{}_{1}F_{1}\left(  a;b;x\right)  \right]  ^{\prime}=\left(
b-1\right)  x^{b-2}{}_{1}F_{1}\left(  a;b-1;x\right)  \label{EqnA29}%
\end{equation}%
\begin{equation}
\left[  x^{b-a}{}e^{-x}{}_{1}F_{1}\left(  a;b;x\right)  \right]  ^{\prime
}=\left(  b-a\right)  x^{b-a-1}e^{-x}{}_{1}F_{1}\left(  a-1;b;x\right)
\label{EqnA30}%
\end{equation}%
\begin{equation}
\left[  e^{-x}{}_{1}F_{1}\left(  a;b;x\right)  \right]  ^{\prime}=\frac
{a-b}{b}e^{-x}{}_{1}F_{1}\left(  a;b+1;x\right)  \label{EqnA31}%
\end{equation}%
\begin{equation}
\left[  x^{b-1}e^{-x}{}_{1}F_{1}\left(  a;b;x\right)  \right]  ^{\prime
}=\left(  b-1\right)  x^{b-2}e^{-x}{}_{1}F_{1}\left(  a-1;b-1;x\right)
\label{EqnA32}%
\end{equation}%
\[
\left[  e^{-x}\left(  x+a-1\right)  ^{b-a}{}_{1}F_{1}\left(  a;b;x\right)
\right]  ^{\prime}=
\]%
\begin{equation}
\frac{\left(  b-a\right)  \left(  b-a+1\right)  }{b}e^{-x}\left(
x+a-1\right)  ^{b-a-1}{}_{1}F_{1}\left(  a-1;b+1;x\right)  \label{EqnA33}%
\end{equation}%
\begin{equation}
\left[  x^{b-1}\left(  x+a\right)  ^{a-b+1}{}_{1}F_{1}\left(  a;b;x\right)
\right]  ^{\prime}=a\left(  b-1\right)  x^{b-2}\left(  x+a\right)  ^{a-b}%
{}_{1}F_{1}\left(  a+1;b-1;x\right)  \label{EqnA34}%
\end{equation}

Inspection of the equations above for Tricomi and Kummer functions shows that
if we express the left-hand side of the recurrences in the form $\left[
u_{i}\left(  x\right)  U\left(  a;b;x\right)  \right]  ^{\prime}$or $\left[
u_{i}\left(  x\right)  {}_{1}F_{1}\left(  a;b;x\right)  \right]  $, the
functions $u_{i}\left(  x\right)  $ are the same for both Tricomi and Kummer
functions. However, the right-hand sides of all of these recurrences are
different. As the $u_{i}\left(  x\right)  $ are the same, the differential
equations needed for the involution are the same. For $u_{1}\left(  x\right)
\equiv1$, the differential equation is (\ref{EqnA1}) and the other 7
differential equations were given in [3] and are repeated below, as they are
needed for the involution.%
\begin{equation}
u_{2}^{\prime\prime}\left(  x\right)  +\left(  \frac{b-2a}{x}-1\right)
u_{2}^{\prime}\left(  x\right)  +\frac{a\left(  a-b+1\right)  }{x^{2}}%
u_{2}\left(  x\right)  =0\label{EqnA35}%
\end{equation}%
\begin{equation}
u_{3}^{\prime\prime}\left(  x\right)  +\left(  \frac{2-b}{x}-1\right)
u_{3}^{\prime}\left(  x\right)  +\frac{\left(  b-a-1\right)  }{x}u_{3}\left(
x\right)  =0\label{EqnA36}%
\end{equation}%
\begin{equation}
u_{4}^{\prime\prime}\left(  x\right)  +\left(  \frac{2a-b}{x}+1\right)
u_{4}^{\prime}\left(  x\right)  +\frac{\left(  a-1\right)  \left(  a-b\right)
}{x^{2}}u_{4}\left(  x\right)  =0\label{EqnA37}%
\end{equation}%
\begin{equation}
u_{5}^{\prime\prime}\left(  x\right)  +\left(  \frac{b}{x}+1\right)
u_{5}^{\prime}\left(  x\right)  +\frac{b-a}{x}u_{5}\left(  x\right)
=0\label{EqnA38}%
\end{equation}%
\begin{equation}
u_{6}^{\prime\prime}\left(  x\right)  +\left(  \frac{2-b}{x}+1\right)
u_{6}^{\prime}\left(  x\right)  -\frac{a-1}{x}u_{6}\left(  x\right)
=0\label{EqnA39}%
\end{equation}%
\begin{equation}
u_{7}^{\prime\prime}\left(  x\right)  +\left(  \frac{b}{x}+1+\frac{2\left(
a-b\right)  }{x+a-1}\right)  u_{7}^{\prime}\left(  x\right)  +\frac{\left(
a-1\right)  \left(  b-a\right)  \left(  a-b-1\right)  }{x\left(  x+a-1\right)
^{2}}u_{7}\left(  x\right)  =0\label{EqnA40}%
\end{equation}%
\begin{equation}
u_{8}^{\prime\prime}\left(  x\right)  +\left(  \frac{b}{x}-1-\frac{2a\left(
x+b-1\right)  }{\left(  x+a\right)  x}\right)  u_{8}^{\prime}\left(  x\right)
+\frac{a\left(  a-b+1\right)  \left(  b-a-2\right)  }{x\left(  x+a\right)
^{2}}u_{8}\left(  x\right)  =0\text{.}\label{EqnA41}%
\end{equation}

\section{Uniqueness of the Involution for Tricomi and Kummer functions}

The right-hand sides of the integrated recurrences are different for the two
functions, so the uniqueness proof given above for Bessel and Associated
Legendre functions can not be applied directly for Tricomi and Kummer
functions. However, there is a very nice way around this difficulty which also
illustrates the importance of conjugate pairs for the overall structure of
this family of functions. It is possible to scale the Tricomi functions with
Gamma function scaling factors $s\left(  a,b\right)  $ such that the Tricomi
functions obey the same recurrences as the Kummer functions. We define a new
set of functions $V_{i}\left(  a;b;x\right)  \equiv s_{i}\left(  a,b\right)
U\left(  a;b;x\right)  $ and since $s_{i}\left(  a,b\right)  $ is not a
function of $x$, the $V\left(  a;b;x\right)  $ obey the various differential
equations given above. The scalings can be obtained by inspection and the
scaled functions are given below.%
\begin{equation}
V_{1}\left(  a;b;x\right)  =\left(  -1\right)  ^{b}\Gamma\left(  b\right)
U\left(  a;b;x\right)  \text{ for Equation (\ref{EqnA19})}\label{EqnB1}%
\end{equation}%
\begin{equation}
V_{2}\left(  a;b;x\right)  =\Gamma\left(  a-b+1\right)  U\left(  a;b;x\right)
\text{ for Equation (\ref{EqnA20})}\label{EqnB2}%
\end{equation}%
\begin{equation}
V_{3}\left(  a;b;x\right)  =\left(  -1\right)  ^{b}\Gamma\left(  a-b+1\right)
\Gamma\left(  b\right)  U\left(  a;b;x\right)  \text{ for Equation
(\ref{EqnA21})}\label{EqnB3}%
\end{equation}%
\begin{equation}
V_{4}\left(  a;b;x\right)  =\Gamma\left(  a-b+1\right)  U\left(  a;b;x\right)
\text{ for Equation (\ref{EqnA22})}\label{EqnB4}%
\end{equation}%
\begin{equation}
V_{5}\left(  a;b;x\right)  =\left(  -1\right)  ^{b}\Gamma\left(  a-b+1\right)
\Gamma\left(  b\right)  U\left(  a;b;x\right)  \text{ for Equation
(\ref{EqnA23})}\label{EqnB5}%
\end{equation}%
\begin{equation}
V_{6}\left(  a;b;x\right)  =\left(  -1\right)  ^{b}\Gamma\left(  b\right)
U\left(  a;b;x\right)  \text{ for Equation (\ref{EqnA24})}\label{EqnB6}%
\end{equation}%
\begin{equation}
V_{7}\left(  a;b;x\right)  =\frac{\left(  -1\right)  ^{b}\Gamma\left(
b\right)  }{\Gamma\left(  b-a\right)  }U\left(  a;b;x\right)  \text{ for
Equation (\ref{EqnA25})}\label{EqnB7}%
\end{equation}%
\begin{equation}
V_{8}\left(  a;b;x\right)  =\frac{\left(  -1\right)  ^{b}\Gamma\left(
b\right)  }{\Gamma\left(  b-a\right)  }U\left(  a;b;x\right)  \text{ for
Equation (\ref{EqnA26})}\label{EqnB8}%
\end{equation}
It has been verified using Mathematica that the $V_{i}\left(  a;b;x\right)  $
obey the same recurrences as the Kummer functions. Inspection of the scalings
above shows that the scalings match for each conjugate pair but not globally
for all the functions. Hence uniqueness can be proven for each conjugate pair
using Wronskians in the same manner as for Bessel functions shown above. Since
the involution carries across its own scaling to its conjugate function, there
would be a contradiction if the scalings did not match, but they do, and all
is well.

\section{Indefinite integrals for Tricomi functions}

The integrated recurrence relations given by Equations (\ref{EqnA19}) to
(\ref{EqnA26}) can be integrated immediately to give the indefinite integrals:%
\begin{equation}%
{\displaystyle\int}
U\left(  a;b;x\right)  dx=-\frac{U\left(  a-1;b-1;x\right)  }{a-1}%
\label{EqnC1}%
\end{equation}%
\begin{equation}%
{\displaystyle\int}
x^{a-2}U\left(  a;b;x\right)  dx=\frac{x^{a-1}U\left(  a-1;b,x\right)
}{\left(  a-1\right)  \left(  a-b\right)  }\label{EqnC2}%
\end{equation}%
\begin{equation}%
{\displaystyle\int}
x^{b-1}U\left(  a;b;x\right)  dx=\frac{x^{b}}{b-a}U\left(  a;b+1;x\right)
\label{EqnC3}%
\end{equation}%
\begin{equation}%
{\displaystyle\int}
x^{b-a-2}e^{-x}U\left(  a;b;x\right)  dx=-x^{b-a-1}e^{-x}U\left(
a+1;b;x\right) \label{EqnC4}%
\end{equation}%
\begin{equation}%
{\displaystyle\int}
e^{-x}U\left(  a;b;x\right)  dx=-e^{-x}U\left(  a;b-1;x\right) \label{EqnC5}%
\end{equation}%
\begin{equation}%
{\displaystyle\int}
x^{b-1}e^{-x}U\left(  a;b;x\right)  dx=-x^{b}e^{-x}U\left(  a+1;b+1;x\right)
\label{EqnC6}%
\end{equation}%
\begin{equation}%
{\displaystyle\int}
\left(  a+x\right)  ^{b-a-3}e^{-x}U\left(  a;b;x\right)  dx=-\left(
a+x\right)  ^{b-a-2}e^{-x}U\left(  a+1;b-1;x\right) \label{EqnC7}%
\end{equation}%
\begin{equation}%
{\displaystyle\int}
x^{b-1}\left(  a-1+x\right)  ^{a-b-2}U\left(  a;b;x\right)  dx=\frac
{x^{b}\left(  a-1+x\right)  ^{a-b-1}}{\left(  a-1\right)  (b-a)\left(
a-b-1\right)  }U\left(  a-1;b+1;x\right)  \text{.}\label{EqnC8}%
\end{equation}
These integrals constitute a fundamental set for the Tricomi functions.
Equations (\ref{EqnC1})-(\ref{EqnC3}) are given in [6]

\section{Comments and conclusions}

Two new integrated recurrence relations have been given for Tricomi functions,
completing the set eight relations linking the Tricomi function ${}U\left(
a;b;x\right)  $ to the eight contiguous functions ${}U\left(  a\pm
1;b\pm1;x\right)  ,{}U\left(  a\pm1;b;x\right)  $ and ${}U\left(
a;b\pm1;x\right)  $. These results were obtained directly from some three-term
contiguous relations given in [6]. Similar results were given recently for the
Kummer Function $_{1}F_{1}\left(  a;b;x\right)  $, with both functions being
independent solutions of Kummer's equation [5]. The involution transformation
presented in [3] has been examined in much more detail and has been applied
here to Cylinder functions, Associated Legendre functions and Tricomi
functions. For all functions to which it applies, it splits the integrated
recurrences into conjugate pairs, with the parameter shifts for conjugate
pairs being reversed. The uniqueness of the involution is easily proven for
Cylinder functions and Associated Legendre functions using Wronskians. The
uniqueness proof is more intricate and also more interesting for Kummer and
Tricomi functions because although both functions obey Kummer's equation, they
have different recurrences. This detail was overcome by using Gamma functions
to scale the Tricomi functions to match the Kummer recurrences. It was found
that the scaling was the same for each conjugate pair, but global scaling is
not possible. As the involution always exists, it will necessarily transfer
the scaling between conjugate pairs and this is indeed the case. This means
that the involution is an important symmetry of the solutions of the various
differential equations to which it applies.

If equations with three parameters are considered, such as the Gauss
hypergeometric function, the grid of parameter shifts is three-dimensional and
it becomes clear that the involution is an inversion through the central point
and not a reflection. This case will be considered in a future paper, as will
Whittaker functions, which can be directly mapped from the Kummer and Tricomi functions.

The uniqueness of the involution means that each contiguous function is a
single occupancy slot, with only one integrated solution for each of the
independent solutions of the differential equation. This is structurally
reminiscent of the Pauli Exclusion Principle in Quantum Mechanics. All new
equations have been checked using Mathematica [9]

\section{Disclosure statement}

No potential conflict of interest was reported by the author.

\end{document}